\documentclass{amsart}
\usepackage{amssymb}   
\usepackage{amsmath}
\usepackage{amsthm}

\def\QQ{{\mathbb Q}}
\def\ZZ{{\mathbb Z}}
\def\Snm{\mathfrak S_n^{(m)}} 
\def\v{\varepsilon}  
\def\={\;=\;} 
\def\O{\text O} 
\def\C{\text{\bf C}}

\newtheorem{theorem}{Theorem}[section]
\newtheorem{lemma}[theorem]{Lemma}
\newtheorem{proposition}[theorem]{Proposition}
\newtheorem{corollary}[theorem]{Corollary}

\newtheorem{definition-lemma}[theorem]{Definition-Lemma}

\theoremstyle{definition}

\newtheorem{example}[theorem]{\bf Example}

\theoremstyle{remark}
\newtheorem{remark}[theorem]{\bf Remark}

\begin{document}

\title[Cycle Relations on Jacobian Varieties]
{Cycle Relations on Jacobian Varieties}
\author{Gerard van der Geer $\,$}
\address{Korteweg-de Vries Instituut, Universiteit van
Amsterdam, Plantage Muidergracht 24, 1018 TV Amsterdam, The Netherlands}
\email{geer@science.uva.nl}
\author{$\,$ Alexis Kouvidakis}
\address{Department of Mathematics, University of Crete, GR-71409 Heraklion, Greece}
\email{kouvid@math.uoc.gr}

\begin{abstract}
By using the Grothendieck-Riemann-Roch theorem we derive cycle relations modulo
algebraic equivalence in the Jacobian of a curve. The relations generalize
the relations found by Colombo and van Geemen and are analogous to but simpler
than the relations recently found by Herbaut. In an appendix due to Zagier it
is shown that these sets of relations are equivalent.
\end{abstract}

\maketitle

\begin{section}{Introduction}\label{sec:intro}
Beauville showed in \cite{B1} that the Chow ring with rational coefficients
of an abelian variety 
possesses a double grading 
$CH_{\QQ}(X) =\oplus \, CH^i_{(j)}(X)$ where $i$ refers to the codimension
and $j$ refers to the action of the integers: 
$CH_{(j)}^i(X)=\{ x \in CH^i(X): k^*(x)=k^{2i-j}\, x\}$. The quotient $A(X)$
of the Chow ring modulo algebraic equivalence inherits this double grading
$A(X)=\oplus \, A^i_{(j)}(X)$ 
and carries two multiplication laws, the intersection product $x \cdot y$ and 
the Pontryagin product $x \ast y$. 

If $X={\rm Jac}(C)$ is the
Jacobian of a curve $C$ of genus $g$ then we can decompose the class $[C]$ of the
image of the Abel-Jacobi map of $C$ as $[C]= \sum_{j=0}^{g-1} C_{(j)}$
with $C_{(j)} \in A^{g-1}_{(j)}({\rm Jac}(C))$.  
Colombo and van Geemen proved (cf.\ \cite{C-vG}) that for a curve $C$ with
a map of degree $d$ to
${\bf P}^1$ the component $C_{(j)}$ vanishes for $j\geq d-1$. 
In \cite{Herb} Herbaut extended this result and 
found cycle relations for curves having a $g_d^r$, i.e., a
linear system of degree $d$ and projective dimension $r$, with $r\geq 2$.

It is the purpose of this note to show that one can use the Grothendieck-Riemann-Roch 
theorem to derive in an easy way the Colombo-van Geemen result as well
as simple relations of higher degree.

Let $C$ be a smooth projective curve of genus $g$ over an algebraically closed
field~$K$. 

\begin{theorem}\label{MainThm} 
If $C$ has a base point free $g^r_d$ then
\[
\sum_{a_1+ \cdots +a_r=N} (a_1+1)! \cdots (a_r+1)! \;
C_{(a_1)}\ast \cdots \ast C_{(a_r)}=0,
\]
for every $N \geq d-2r+1$.
\end{theorem}
For $N=d-2r+1$ this relation coincides with Herbaut's relation (cf., \cite{Herb}, Thm 1 and Thm 8). 
For higher values of $N$ they are in general different, 
but in an appendix we present a proof by Don Zagier that shows that these
sets of relations are equivalent.

\end{section}

\begin{section}{Preliminaries}
Let $C$ be a smooth projective curve of genus $g$ over an algebraically closed
field $K$. We suppose that the curve $C$ has a base-point free
linear system $g_d^r$ of degree $d$ and projective dimension $r$.
This defines a morphism $\gamma: C \to {\bf P}^r$.
Let $J={\rm Jac}(C)$ be the Jacobian of $C$.

We consider the incidence variety $Y \subset C \times \hat{\bf P}^r$ defined by
$$
Y=\{ (p,\eta) \in C \times \hat{\bf P}^r \colon \gamma(p) \in \eta \},
$$
where $\hat{\bf P}^r$ is the dual projective space of ${\bf P}^r$.
It has dimension $r$ and possesses the two projections $\tilde{\phi}$
and $\tilde{\alpha}$ onto $C$ and $\hat{\bf P}^r$. Note that
$\tilde{\alpha}$ is finite of degree $d$ and
$\tilde{\phi}$ is a ${\bf P}^{r-1}$-fibration.
We shall write ${\bf P}^r$
for $\hat{\bf P}^r$.

We have the following diagram of morphisms
$$
\begin{matrix}
& {\bf P}^r &{\buildrel v \over \longleftarrow} & {\bf P}^r\times J &
{\buildrel \alpha \over \longleftarrow}  & Y\times J&
{\buildrel \pi \over \longrightarrow} & Y & {\buildrel \tilde{\alpha} \over \to}
 & {\bf P}^r  \cr
&&&\downarrow{p} && \downarrow{\phi} && \downarrow{\tilde{\phi}} \cr
&&&J & {\buildrel q \over \longleftarrow} &
C \times J & {\buildrel \tilde{\pi} \over \longrightarrow} & C \, ,
\cr
\end{matrix}
$$
where the morphisms $v, p, q, \tilde{\pi}$ and $\pi$ are projections
and $\alpha= \tilde{\alpha}\times {\rm id}_J$ and 
$\phi=\tilde{\phi} \times {\rm id}_J$.

Let $P$ be the Poincar\'e bundle on $C \times J$ and
set $L:=\phi^* P$, a line bundle on $Y \times J$. Put $\ell:=c_1(L)$
and $\Pi:=c_1(P)$.

The Chow ring of ${\bf P}^r \times J$ is generated over $CH^*(J)$
by the class $\xi=v^*{h}$ with $h$ a hyperplane in ${\bf P}^r$ with
$\xi^{r+1}=0$.
For a class
$\beta \in CH^*({\bf P}^r \times J)$ we have the relation (cf.\ \cite{Fulton},
(Thm.\ 3.3, p.\ 64))
$$
\beta= \sum_{i=0}^r \beta_i \, \xi^{r-i} \qquad {\rm with} \quad
\beta_i=p_*(\beta  \cdot \xi^i)\, ,\eqno(1)
$$
where by abuse of notation we write here and hereafter 
$\beta_i$ for $p^*(\beta_i)$.

We let
$x=\alpha^*(\xi)$ be the pull back of $\xi$. We let
$\rho=\pi^* \tilde{\phi}^*{\rm (point)}$
be the pull back class of a point on $C$. We work in the
 Chow ring up to algebraic equivalence. There we have the relations
$$
x^r=d \, \pi^*({\rm point}), \quad x^{r+1}=0, \quad \rho^2=0, \quad
x^{r-1}\rho= \pi ^* ({\rm point}).
$$

Recall the Fourier transform $F: A(X) \to A(X)$ for a principally polarized
abelian variety $(X,\theta)$ of dimension $g$, cf.\ \cite{B1, B2}. 
It has the properties i) $F \circ F= (-1)^g (-1)^*$, 
ii) $F(x \ast y)= F(x) \cdot F(y)$
and $F(x \cdot y)= (-1)^g F(x)\ast F(y)$,
iii) $F(A^i_{(j)}(X))= A^{g-i+j}_{(j)}(X)$.

We have the relation $q_*(e^{\Pi})=F[C]$, cf.\ \cite{B2}, Section 2.
Comparing terms gives that $F[C_{(j)}]= (1/(j+2)!) \, q_*(\Pi ^{j+2})$ 
for $j=0, \ldots ,g-1$.
Note, also, that $q_*1=q_*\Pi=0$. More generally, extending scalars to $\QQ$ we 
have the relation
$$
q_*(e^{k\Pi})=k^{2g}F[(k^{-1})^*C] \qquad \hbox{\rm for $k\in \ZZ_{\geq 1}$}.
$$
In fact, writing $[C]=\sum_{j=0}^{g-1} C_{(j)}$ we have
$(k^{-1})^*[C]=\sum_{j} k^{j+2-2g} C_{(j)}$, hence
$$
k^{2g}F[(k^{-1})^*C]= F[\sum_{j} k^{j+2}C_{(j)}]=
\sum_j k^{j+2} \, q_*(\Pi^{j+2}/(j+2)!)=q_*(e^{k\Pi}).
$$
\end{section}

\begin{section}{The Proof}

We shall prove that
if $C$ has a base point free $g^r_d$ then
\[
\sum_{a_1+ \cdots +a_r=N} (a_1+1)! \cdots (a_r+1)! \;
F[C_{(a_1)}]\cdots F[C_{(a_r)}] =0,\eqno(2)
\]
for every $N \geq d-2r+1$.

We are going to apply the Grothendieck-Riemann-Roch theorem to the morphism
$\alpha$ and the line bundle $L$. For $k\geq 1$ we put $V_k:=\alpha_*(L^{\otimes k})$. 
Since $\alpha$
is a finite morphism of degree $d$ this is a vector bundle of rank $d$ and we get
$$
{\rm ch}(V_k)={\rm ch}(\alpha_{!}L^{\otimes k})=\alpha_*(e^{k \ell}{\rm td}_{\alpha}),
$$
with ${\rm td}_{\alpha}$ the Todd class of the morphism $\alpha$.

The Todd class ${\rm td}_{{\alpha}}$ is algebraically
equivalent to a class of the form $A(x)+B(x)\, \rho$. 
Here $A=\sum_{j=0}^{r-1} a_j x^j$ and $B=\sum_{j=0}^{r-1}b_j x^j$ are polynomials
in $x$ and $a_0=1$. In fact,
${\rm td}_{\alpha}$ is the pull back under $\pi$ of ${\rm td}_{\tilde{\alpha}}$,
an element of $A(Y)$. The ring $A(Y)$ is generated as an $A(C)$-module by
$1,x_1,\ldots,x_1^{r-1}$ with $x_1=\tilde{\alpha}^*(h)$ and $A(C)$ is generated by $1$ and 
the class of a point.
 
\begin{proposition}\label{chernprop} 
For $k\in {\ZZ}_{\geq 1}$ we have in $A({\bf P}^r \times J)$ the relation
$$
{\rm ch}(V_k) = d A(\xi ) +\xi B(\xi) + k^{2g}F[(k^{-1})^*C]  \xi A(\xi)\,   .
$$
In particular, all ${\rm ch}_j(V_k)$ are divisible by $\xi$ for $j\geq 1$.
\end{proposition}

Before we give the proof of Proposition  \ref{chernprop} we state a 
corollary and a lemma.
\smallskip

Proposition \ref{chernprop}  gives an expression of the Chern characters of 
the bundles $V_k$. We can express the Chern classes of the bundles $V_k$ 
of rank  $d$ in terms of the Chern characters by using the well-known
formula (cf.\ \cite{MacDonald}, ${\rm ch.}\, {\rm I}\; (2.10')$)
\[
1+c_1(V_k)\, t + \cdots + c_d(V_k)\, t^d= 
{\rm exp} ( \sum _{j\geq 1} (-1)^{j-1}(j-1)!\; 
{\rm ch}_j(V_k) \; t^j ). \eqno(3)
\]
Formula (3) combined with Proposition \ref{chernprop}   will 
give us the vanishing relations we are asking for. For example, applying these formulas
for $r=1$ and $k=1$ immediately gives us the Theorem of 
Colombo-van Geemen \cite{C-vG} as we now show.

\begin{corollary} 
If $C$ has a $g^1_d$ then $C_{(j)}=0$ for $j\geq d-1$.
\end{corollary}
\begin{proof} Put $V=V_1$. 
We see ${\rm ch}(V)=d+n\xi+F[C]\xi$ for some $n$ (actually $n=1-d-g$). Since  
${\rm ch}_j(V)$ is  divisible by $\xi $ for $j\geq 1$ and $\xi ^2=0$, formula 
(3) becomes in this case:  $1+c_1(V)\, t + \cdots + c_d(V)\, t^d =
1 + {\rm ch}_1(V)\, t - {\rm ch}_2(V)\, t^2+ \cdots + (-1)^{j-1}(j-1)!\,  
{\rm ch}_j(V) \; t^j+ \cdots $. Therefore ${\rm ch}_j(V)$ vanishes for $j>d$.
Hence $F[C]$ has no terms of codimension $\geq d$.  Since $F[C_{(j)}]$ is of 
codimension $j+1$, it follows that $C_{(j)}=0$ for all $j\geq d-1$.
\end{proof}

\begin{lemma}
In $A({\bf P}^r \times J)$ the following relations hold for
$\nu \geq 0$ :
$$
\alpha_*(\ell^{\mu} \cdot x^{\nu})=\begin{cases}
 q_*(\Pi^{\mu}) \, \xi^{\nu+1} & \mu>0, \\
d \, \xi^{\nu} & \mu=0 \\
\end{cases}
\;\;\; \mbox{ and } \;\;\;
\alpha_*(\ell^{\mu} \cdot x^{\nu} \cdot \rho )=\begin{cases} 0  & \mu>0, \\
\xi^{\nu +1}& \mu=0\, . \\
\end{cases}
$$
\end{lemma}

\begin{proof}
For the first relation: By (1) the coefficient of $\xi^{r-j}$ is given by
\begin{eqnarray*}
p_*(\alpha_* (\ell ^{\mu} \, x^{\nu}) \xi^j) &=& p_*(\alpha_*(\ell^{\mu }
x^{\nu} \alpha^*(\xi)^{j})) = p_*(\alpha_*(\phi^*(\Pi)^{\mu} x^{\nu+j} ))\\
&=& q_*(\phi_*(\phi^*(\Pi)^{\mu} x^{\nu+j} )) = q_*(\Pi^{\mu} \phi_*(x^{\nu+j})).
\end{eqnarray*}
If $\nu+j=r$ then $x^r$ is algebraically equivalent to $d$ times ${\rm point}\times J$
and since $\Pi$ is algebraically equivalent to $0$ on ${\rm point}\times J$
we get that any term
with $\mu>0$ and $\nu+j = r$ contributes $0$. The term with $\nu+j=r-1$
contributes $q_*(\Pi^{\mu} \phi_*(x^{r-1}))=q_*(\Pi^{\mu})$ since
$\phi_*(x^{r-1}))=1_{C \times J}$. If $\nu+j<r-1$ or $\nu+j>r$ then we get
$q_*(\Pi^{\mu} \phi_*(x^{\nu+j}))=q_*(\Pi^{\mu} \, 0)=0$. Finally, if $\mu=0$ we use
that $x=\alpha^*\xi$, hence $\alpha_*(x^{\nu})=d\, \xi^{\nu}$.

For the second relation: Observe that $\ell ^{\mu } \cdot \rho =0$, if $\mu \geq 1$. 
Indeed, $\ell = \phi ^* \Pi$, $\rho = \phi ^* (  {\rm point}  \times J )$ and 
$\Pi \cdot ( {\rm point}  \times J) =0$. When $\mu =0$ we have
$\alpha _* (x^{\nu} \rho)= \xi ^{\nu} \alpha _*\rho = \xi^{\nu +1}$.

\end{proof}

We now give the proof of Proposition  \ref{chernprop}.

\begin{proof}
We have $\alpha_*(\ell^{\mu} x^{\nu} \rho)=0 \quad \hbox{\rm for all $\mu\geq 1$}$, 
$\nu \geq 0$. So in the contributions $\alpha_*(e^{k\ell} {\rm td_{\alpha}})$ 
we get contributions of the form $\alpha_*(e^{k\ell} A(x))$ and
$\alpha_*(\rho B(x))$ only. We have
\begin{eqnarray*}
\alpha_*(e^{k\ell} A(x)) &=& \alpha_*(A(x)) + \sum _{\mu \geq 1} 
\frac{k^{\mu}}{\mu !} \alpha_* (\ell ^{\mu} A(x)) 
 =    dA(\xi) + \sum _{\mu \geq 1}  \frac{k^{\mu}}{\mu !} q_* \Pi^{\mu} \xi A(\xi) \\
&=&   dA(\xi) + q_*(e^{k\Pi})  \xi A(\xi) =dA(\xi) + k^{2g}F[(k^{-1})^*C]  \xi A(\xi).   
\end{eqnarray*}
On the other hand, $\alpha_*(\rho B(x))=\xi B(\xi)$. 
\end{proof}
By using the relation 
$$k^{2g}F[(k^{-1})^*C]  = q_*(e^{k\Pi})  = \sum _{\mu \geq 0} \frac{k^{\mu+2}}{\mu+2 !} 
q_* \Pi^{\mu+2} = \sum _{\mu \geq 0}  k^{\mu+2} \, F[C_{(\mu)}],
$$ 
we get the following corollary of Proposition  \ref{chernprop}.
\begin{corollary}
We put $a_j= b_j =0$ for every $j \geq r$. With $A=\sum _{j=0}^{r-1}a_jx^j$ 
where $a_0=1$ and 
$B=\sum _{j=0}^{r-1}b_jx^j$ we have for $j \geq 1,\, k \geq 1$ that
\[
{\rm ch}_j(V_k)= (d\,a_j+b_{j-1})\, \xi^j 
+ \sum_{m=1}^{j-1} a_{m-1}\, k^{j-m+1} \, F[C_{(j-m-1)}] \; \xi ^m \, .
\]
\end{corollary}

With $j \geq 1$ we  put
\[ 
{\rm ch}_j(V_{k}) =  A_1(j) \xi ^1+ \cdots + A_{r}(j) \xi^{r},
\]
where $A_m(j)$ is of codimension $j-m$. Please note that ${\rm ch}_j(V_k)$ 
is divisible by $\xi $ for $j \geq 1$ by Prop.\ \ref{chernprop}.
\begin{remark}
The coefficient $A_m(j)$ depends on $k$, but for simplicity of notation
we do not involve the index $k$ in the notation.
\end{remark}
 Then, for every $j \geq 1$ and $1 \leq m \leq r$, we have
$$
A_m(j) = \begin{cases}
  d\,a_j+b_{j-1}   &  m=j, \\
 a_{m-1}\, k^{j-m+1}F[C_{(j-m-1)}] &   m<j, \\  0   & m>j. 
\end{cases}
$$

Since ${\rm rank}(V_k)=d$ the coefficient of $t^{M+1}$ of the 
right hand side of (3) must be zero for every $M \geq d$. Let us write
\[ 
F(t)= \sum _{j\geq 1} (-1)^{j-1}(j-1)!\; {\rm ch}_j(V_k) \; t^j \;. 
\]
Note that $F(t)^i=0$ for every $i \geq r+1$ because ${\rm ch}_j(V_k)$ 
is divisible by $\xi$ for $j \geq 1$. Therefore the right hand side of 
(3) is equal to
$ \sum_{i=0}^r (1/i!) \,  F(t)^i$.
The coefficient of $t^{M+1}$ in the polynomial $F(t)^i$, for $i \geq 1$, is equal to
$$
(-1)^{M+1-i} \sum_{\alpha _1+ \cdots +\alpha _i=M+1} 
(\alpha _1-1)!\cdots  (\alpha _i-1)!\;
{\rm ch}_{\alpha _1}(V_k) \cdots {\rm ch}_{\alpha _i}(V_k)\, .  
$$
We denote the expression 
$(\alpha _1-1)!\cdots  (\alpha _i-1)!$ by $\alpha \{1,i\}$.
Then we have that
\[
\sum_{i=1}^{r} \frac{(-1)^{i}}{i!} \sum_{\alpha _1+ \cdots +\alpha _i=M+1}
 \alpha \{1,i\}\; {\rm ch}_{\alpha _1}(V_k) \cdots  {\rm ch}_{\alpha _i} (V_k) = 0\,    
\]
for $M \geq d$ and so, 
\[
\sum_{i=1}^r \frac{(-1)^{i}}{i!} 
\sum_{\alpha _1+ \cdots +\alpha _i=M+1} \alpha \{1,i\}\;  
[\sum_{m_1=1}^r A_{m_1}(\alpha_1) \xi ^{m_1}] \cdots
[\sum_{m_i=1}^rA_{m_i}(\alpha_i)\xi ^{m_i}]=0 \, . 
\]
Now the LHS of this is easily seen to be equal to
\[
\sum_{m=1}^r \sum_{i=1}^m  \frac{(-1)^{i}}{i!} 
\sum_{\alpha _1+ \cdots +\alpha _i=M+1} 
\alpha \{1,i \}\; \sum_{m_1+\cdots +m_i=m }   
[A_{m_1}(\alpha_1) \cdots A_{m_i}(\alpha_i)] \xi ^m\, . 
\]
We therefore have  for every $m=1, \ldots , r$ and $M\geq d$ that
\[
\sum_{i=1}^m \frac{(-1)^{i}}{i!} \sum_{m_1+\cdots +m_i=m \;}   
\sum_{\,\alpha _1+ \cdots +\alpha _i=M+1} \alpha \{1,i\}\;  
A_{m_1}(\alpha_1) \cdots A_{m_i}(\alpha_i) =0\, .
\]
With $M \geq d$ the case $m=r$ gives the relation
\[
\sum_{i=1}^r \frac{(-1)^{i}}{i!} \sum_{m_1+\cdots +m_i=r \;}  
 \sum_{\,\alpha _1+ \cdots +\alpha _i=M+1} \alpha \{1,i\}\;
A_{m_1}(\alpha_1) \cdots A_{m_i}(\alpha_i) =0\, .
\]
If we write 
\[ B_M(i)=  \sum_{m_1+\cdots +m_i=r \;}   
\sum_{\,\alpha _1+ \cdots +\alpha _i=M+1} \alpha \{1,i\}\;
A_{m_1}(\alpha_1) \cdots A_{m_i}(\alpha_i) 
\]
then this relation becomes $\sum_{i=1}^r ((-1)^{i}/{i!}) B_M(i)  =0$ 
for every $M \geq d$.  We analyze the dependence on $k$. 
\begin{proposition}
We write $\sum _{i=1}^r \frac{(-1)^{i}}{i!} B_M(i) = 
\sum _s \Gamma_s k^s$ as a polynomial in $k$. With $M \geq d$ we have that 
$\Gamma_s=0$ for $s > M+1$ and
\[\Gamma_{M+1}=  \frac{(-1)^r}{r!} \sum_{\alpha _1+ \cdots +\alpha _r=M-2r+1} 
(\alpha_1+1)! \cdots (\alpha _r+1)! \; F[C_{a_1}]\cdots F[C_{a_r}] \, . 
\]
\end{proposition}

\noindent
\begin{proof}
If $A_{m_1}(\alpha_1) \cdots A_{m_i}(\alpha_i)$ contains a factor with $m_j > \alpha _j$ 
then it vanishes. Otherwise, since $A_{m_j}(\alpha _j) =  a_{m_j-1}\, 
k^{\alpha _j-m_j+1}F[C_{\alpha_j-m_j-1}]$, except when $m_j=\alpha _j$, in which case 
$A_{\alpha_j}(\alpha _j)= d\, a_{\alpha_j}+b_{\alpha _j-1}$, the power of $k$ contained in
 $A_{m_1}(\alpha_1) \cdots A_{m_i}(\alpha_i)$ is equal to 
$ (\alpha _1 + \cdots +\alpha _i)-(m_1 + \cdots +m_i)+ \nu = M+1-r+\nu$, where 
$\nu $ is given by $\nu =\# \{m_j \neq \alpha j, \; j=1, \ldots , i \}\leq i$.  Now if $i < r$
the above number is $<M+1$.  On the other hand, if $i=r$ then  $m_1=\cdots =m_r=1$ 
(since $m_i\geq 1$) and therefore 
\[
B_M(r)=   \sum_{\,\alpha _1+ \cdots +\alpha _r=M+1} \alpha \{1,r\}\;
A_{1}(\alpha_1) \cdots A_{1}(\alpha_r)\, . 
\]
Again, if a term of the above sum contains a factor $A_1(1)$, then the power of $k$ 
contained in this term is $<M+1$. On the other hand, the sum of the terms with no factor 
of the form $A_1(1)$ is
\[
B'_M(r)= \sum_{a_{\mu} \geq 2,\,\;\alpha _1+ \cdots +\alpha _r=M+1} \alpha \{1,r\}\;
A_{1}(\alpha_1) \cdots A_{1}(\alpha_r) \, .
\]  
Since  $\alpha _\mu \geq 2$ we have that $A_1(\alpha _{\mu })=k^ {\alpha_{\mu }} \; 
F[C_{\alpha_{\mu} -2}] $. 
Therefore
\begin{eqnarray*}
B'_M(r) &=& k^{M+1} \sum_{\;\alpha_{\mu } \geq 2,\,\; \alpha _1+ \cdots +\alpha _{r}=M+1} 
\;\alpha \{1,r\}\; F[C_{\alpha_1 -2}]  \cdots F[C_{\alpha_{r} -2}] \\
 &=& k^{M+1} \sum_{\;\alpha_{\mu } \geq 0,\;\, \alpha _1+ \cdots +\alpha _{r}=M-2r+1} 
\;(\alpha_1+1)! \cdots (\alpha_{r}+1)! \; F[C_{\alpha_1}]  \cdots F[C_{\alpha_{r}}] \, . 
\end{eqnarray*}
\end{proof}

We now prove Theorem \ref{MainThm}
\begin{proof}
Since we are working with ${\mathbb Q}$-coefficients, we conclude 
that the coefficient of $k^{M+1}$ (which is the maximum degree of 
$k$ involved) must be zero for $M \geq d$, which is relation (2) with 
$N=M-2r+1\geq d-2r+1$.  By applying the Fourier transform 
Theorem \ref{MainThm} follows. 
\end{proof}
\end{section}

\begin{section}{Appendix provided by Don Zagier}
\begin{subsection}{A combinatorial identity}
We set
$$ B_d(a_1,\dots,a_r)\=\sum_{i_1,\dots,i_r\ge1}(-1)^{d-i_1-\cdots-i_r}
\binom d{i_1+\cdots+i_r}\,i_1^{a_1}\cdots i_r^{a_r}\, .
$$
Using
$\sum\limits_{d\=0}^\infty(-1)^{d-i}\binom di\,u^d=\dfrac{u^i}{(1+u)^{i+1}}$ 
we find that $ \sum_{d=0}^\infty B_d(a_1,\dots,a_r)\,u^d$ equals
$$
 \sum_{i_1,\dots,i_r\ge1}i_1^{a_1}\cdots i_r^{a_r}\frac{u^{i_1+\cdots+i_r}}{(1+u)^{i_1+
\cdots+i_r+1}} \=\frac1{1+u}\,P_{a_1+1}(u)\cdots P_{a_r+1}(u)\,,
$$
where $P_n(u)$ is the power series
$$ P_n(u)\;:=\;\sum_{i=1}^\infty\,i^{n-1}\,
\biggl(\frac u{1+u}\biggr)^i\quad\in\;\mathbb
 Z[[u]]\qquad(n\ge1)\,. 
$$
\begin{lemma}
\begin{itemize}
\item[(i)] $P_n(u)$ is a polynomial of degree $n$.  More precisely,
we have 
$P_n(u)=\sum_{m=1}^n(m-1)!\,\Snm u^m$ where each $\Snm$ is a positive integer.
\item[(ii)] $P_n(-1)=0$ for $n>1$.
\item[(iii)] For each $n\ge1$ one has the Laurent series identity
$$ 
\frac{(n-1)!}{[\log(1+x)]^n}\=P_n\bigl(\frac{1}{x}\bigr)\,+\,\text{\rm O}(x)
\qquad{\rm in}\quad \mathbb Q[x^{-1},x]]\,.\eqno(4)
$$
\end{itemize}
\end{lemma}
\begin{example}
The first five values of $P_n(u)$ are $u$, $u+u^2$, $u+3u^2+2u^3$,
$u+7u^2+12u^3+6u^4$ and $u+15u^2+50u^3+60u^4+24u^5$.  We have
$$
\frac{4!}{[\log(1+x)]^5}\=\frac{24}{x^5}+\frac{60}{x^4}+\frac{50}{x^3}+
\frac{15}{x^2}+\frac1x+0  -\frac x{252}+\frac{x^2}{504}-\frac{19x^3}{30240}
-\frac{x^4}{20160}+\frac{53x^5}{147840}-\cdots\;.
$$
\end{example}
\begin{proof}
The easiest approach (as usual!) is to use generating functions.  We have
$$ 
\sum_{n=1}^\infty P_n(u)\,\frac{t^{n-1}}{(n-1)!}\=\sum_{i=1}^\infty\biggl(\frac{ue^t}{1+u}
\biggr)^i\=\frac{ue^t}{1\,-\,u(e^t-1)} \eqno(5)
$$
and hence
$$ 
\frac{1}{(n-1)!}\,\C_{u^m}\bigl[P_n(u)\bigr]\= 
\C_{t^{n-1}}\bigl[e^t(e^t-1)^{m-1}\bigr]\qquad\;(m,\,n\ge1)\,. 
$$
(Here $\C_{x^i}[\Phi(x)]$ denotes the coefficient of $x^i$ in a polynomial, power series or Laurent series $\Phi(x)$.)
This clearly vanishes for $n<m$, showing that $P_n$ is a polynomial of 
degree $\le n$, and we also get the explicit formula for $P_n$ given in~(i), with
$$
\Snm\=\frac{n!}{m!}\,\C_{t^n}\bigl[(e^t-1)^m\bigr]
\=\C_{x^{n-m}}\biggl[\prod_{j=1}^m\frac1{1-jx}\biggr]
\=\frac1{m!}\sum_{k=0}^m(-1)^{m-k}\binom mkk^n
$$
being the Stirling number of the second kind 
(= number of ways of partitioning a set of $n$ elements into
$m$ non-empty subsets). The right-hand side of~(5) reduces to $-1$ at $u=-1$, proving~(ii). 
For~(iii), we use the residue theorem and the substitution $e^t=1+x$ to get
$$ 
\frac1{(n-1)!}\,\C_{u^m}\bigl[P_n(u)\bigr]\= 
\text{Res}_{t=0}\biggl[\frac{e^t(e^t-1)^{m-1}\,dt}{t^n}\biggr]
\=  \text{Res}_{x=0}\biggl[\frac{x^{m-1}\,dx}{[\log(1+x)]^n}\biggr]
$$
for $m,\,n\ge1$.
\end{proof}
\end{subsection}
\begin{subsection}{The equivalence of the two sets of relations.}
Let $R$ be the $\mathbb Q$-subalgebra of
the Chow algebra generated (with respect to the Pontryagin product) by the $C_{(j)}$ ($0\le j\le g-1$), bigraded
by $C_{(j)}\in R_{1,j}$, so that $R_{ij}\subseteq A^{g-i}_{(j)}\,$.  We define two polynomials
$$  G(t)\=\sum_{a=0}^{g-1}\,(a+1)!\,C_{(a)}\,t^{a+2}\;\in\;R[t]\,,\qquad
 H(u,t)\=\sum_{a=0}^{g-1}\,P_{a+2}(u)\,C_{(a)}\,t^{a+2}\;\in\;R[u,t]\,. $$
Then the relations obtained from Thm \ref{MainThm} can be written in the form
$$ 
\deg_t\bigl[G(t)^s\bigr]\;\le\;d-r+s\qquad\text{for $1\le s\le r\,$,} \eqno(6)
$$
while, in view of the formulas in the preceding subsection, 
Herbaut's (cf.\ \cite{Herb}, Thm1) can be written in the form
$$ 
\C_{u^{d-r+s}}\bigl[\frac1{1+u}\,H(u,t)^s\bigr]\=
0\quad\text{in $R[t]$}\qquad\text{for $1\le s\le r$.} \eqno(7) 
$$
We also introduce the {\it strengthened Herbaut relations}
$$ 
\deg_u\bigl[H(u,t)^s\bigr]\;\le\;d-r+s\qquad\text{for $1\le s\le r\,$.} \eqno(8) 
$$
Clearly (8) implies (7). (Note that $H(u,t)^s/(1+u)$ is a polynomial by part~(ii) of
the lemma.)  To see that (6) implies (8), we use equation~(4) to obtain
$$  
G\biggl(\frac{t}{\log(1+x)}\biggr)\=H\bigl(\frac{1}{x},\,t)\,-\,\v(x,t)
$$
with $\v(x,t)=\O(x)$ (in fact $\v(x,t)=\text O(xt^2)$) and hence, assuming (6),
\begin{align}
H\bigl(\frac{1}{x},\,t)^si&\=
\sum_{s'=0}^s\binom s{s'}\,G\biggl(\frac{t}{\log(1+x)}\biggr)^{s'}\,\v(x,t)^{s-s'}
\tag{9} 
\\
&\=\sum_{s'=0}^s\O\bigl(\frac{1}{x^{d-r+s'}}\bigr)\,\O\bigl(x^{s-s'}\bigr)\=\O\bigl(
\frac{1}{x^{d-r+s}}\bigr)\notag
\end{align}
as $x\to0$, proving (8). To prove that (7) implies (6), we use induction on~$s$. 
Assume (7) for some $s\le r$ and (6)
for all $s'<s$.  Equation (9) and the inductive assumption give
$$ 
G\biggl(\frac t{\log(1+x)}\biggr)^s\=
H\bigl(\frac1x,\,t)^s\,+\,\O\bigl(\frac1{x^{d-r+s-2}}\bigr)
$$
and hence for $n>d-r+s-2$
$$ 
\C_{t^n}\bigl[G(t)^s\bigr]\cdot\C_{x^{-d+r-s}}
\bigl[\frac x{1+x}\,\frac1{[\log(1+x)]^n}\bigr]
 \= \C_{t^n}\bigl[\C_{u^{d-r+s}}\bigl[\frac1{1+u}\,H(u,t)^s\bigr]\bigr]\=0
$$
by~(7). Differentiating (4) and using part~(i) of the lemma, 
we find that the second factor on the left equals
$\dfrac{(d-r+s)!}{(n-1)!}\,\mathfrak S_{n-1}^{(d-r+s)}$, 
which is non-zero for $n>d-r+s$. Equation (6) follows.

\end{subsection}
\smallskip
\noindent
{\bf Acknowledgements.} Both authors like to thank Don Zagier for providing
a proof of the equivalence of the two sets of relations.
The first author would like to thank the
Department of Mathematics at Heraklion for the hospitality and
the excellent working conditions
during his stay in June/July 2006.
\end{section}

 \end{document}